\documentclass[12pt,a4paper]{amsart}
\usepackage{amssymb,latexsym}
\def\ZZ{\mathbb Z}
\def\RR{\mathbb R}

\def\PP{\mathbb P}

\newtheorem{demo}{D\'emonstration}
\newtheorem{lem}{Lemme}
\newtheorem{theo}[lem]{Th\'eor\`eme}
\newtheorem{cla}[lem] {Affirmation}
 
\newtheorem{prop}[lem]{Proposition} 

\newtheorem{remark}[lem]{Remarque}

\begin{document}
 \title 
[$\lambda _{2g-2}>\frac 14$] 
{Pour  toute surface hyperbolique de genre $g$,
  $\lambda _{2g-2}>\frac 14$}\author{Jean-Pierre Otal 
et Eulalio Rosas}  
\abstract
{We study the influence of the topology of an hyperbolic surface on the number of its Laplace  eigenvalues which are    $\leq \frac 14$. The first result of the paper is  a  ``spectral gap'' statement, namely  its 
title   where $\lambda _j$ is the $j$-th  of the eigenvalues 
$\lambda _0=0<\lambda _1\leq \cdots \leq\lambda _j\leq\cdots $ 
of the Laplace operator and where $g$ is the genus of the surface. 
A classical construction due to Peter Buser shows that this result is sharp.
We give a similar statement for finite volume surfaces.
The methods develop those  of 
\cite {O}, which used  in an essential way the topological  approach
of  Bruno S\'evennec
 to the question of bounding the multiplicity of the second eigenvalue        for  Schr\"odinger operators \cite {Se}}.
 
 \
 
  \noindent{R\'ESUM\'E.} Nous \'etudions l'influence de la topologie d'une surface hyperbolique sur le nombre des valeurs propres de son  Laplacien qui sont $\leq \frac 14$. Le premier r\'esultat de l'article est un \'enonc\'e du type ``trou spectral'', son titre, dans lequel $\lambda _j$ est la $j$-i\`eme des valeurs propres 
  $\lambda _0=0<\lambda _1\leq \cdots \leq\lambda _j\leq\cdots $ 
 du Laplacien et $g$ est le genre de la surface. Une construction  classique d\^ue \`a  Peter Buser montre que ce r\'esultat est optimal. Nous donnons  aussi un \'enonc\'e du m\^eme type pour les surfaces de volume fini. Les m\'ethodes prolongent celles de
  \cite {O}, qui utilisaient de mani\`ere essentielle l'approche topologique par  Bruno S\'evennec de  la question de la majoration de la multiplicit\'e
  de la deuxi\`eme valeur propre  
des op\'erateurs de  Schr\"odinger \cite {Se}.
  \endabstract

\maketitle

\section{\bf Introduction}

Soient $S$ une surface orientable, $\tilde S$ son rev\^etement universel.
On suppose $S$    munie d'une m\'etrique Riemannienne et on note $\Delta$ son Laplacien. Soit 
 $\lambda _0 (\tilde S)$ {\it le bas du spectre du Laplacien} sur   $\tilde S$ ; on  peut
 caract\'eriser 
$\lambda _0 (\tilde S)$ 
comme la borne inf\'erieure des quotients de Rayleigh 
$\displaystyle \frac {\int _{\tilde S}\Vert \nabla \phi \Vert^2}{\int _{\tilde S}\phi ^2}$ lorsque $\displaystyle \phi : \tilde S\to \RR$ d\'ecrit l'ensemble des fonctions $C^\infty$ \`a support compact. 
Quand la m\'etrique est de courbure  constante $-1$, un calcul classique donne\;:  $\displaystyle \lambda _0 (\tilde S)=\frac 14$.

Nous disons qu'une   fonction $C^\infty$ $f: S\to \RR$ est  {\it $\lambda$-propre} si elle v\'erifie : $ \Delta f+\lambda f =0$. 

Lorsque $S$ est compacte, le spectre de l'op\'erateur $-\Delta$  est  un ensemble \;
 discret  $\{0=\lambda _0<\lambda _1\leq \cdots \leq  \lambda _j \leq \cdots \}$ et l'espace $L^2 (S)$ poss\`ede une base orthonorm\'ee de fonctions $\lambda$-propres.

\begin{theo}\label{2g-2}Soit $S$ une surface  compacte
munie d'une m\'etrique analytique de courbure strictement n\'egative.
Alors
la $(2g-2)$-i\`eme valeur propre 
 $\lambda _{2g-2}$     
   est $>\lambda _0 (\tilde S)$.
   \end{theo}  
Pour   une surface  $S$ {\it hyperbolique}, c'est-\`a-dire munie d'une m\'etrique de courbure constante $-1$, Peter Buser a montr\'e  l'in\'egalit\'e $\displaystyle \lambda _{4g-2}> \frac 14$
\cite{Bu1}  \cite{Cha}.  
Buser a aussi montr\'e l'existence d'une constante 
$\epsilon $ 
ind\'ependante de $g$,  de l'ordre de $10^{-12}$ telle que
$\lambda _{2g-2}\geq \epsilon$ \cite {Bu2}.
L'article \cite{Bu1} contient,   
 pour tout $g$,
 des exemples
  de surfaces hyperboliques de genre  $g$
 telles que les  $(2g-3)$ premi\`eres valeurs propres du spectre sont arbitrairement petites\;: ainsi  la conclusion du 
Th\'eor\`eme \ref{2g-2} est optimale quant \`a l'indice de la valeur propre.    Elle est aussi optimale quant \`a la constante $\displaystyle \frac 14$ car 
\cite {Bu1} donne  aussi, pour tout $\epsilon >0$ et pour tout entier $i\geq 0$,  des 
exemples de surfaces  telles que  $\lambda _i $ est inf\'erieur   
\`a $\displaystyle \frac 14+\epsilon$.
 Rappelons maintenant bri\`evement
 comment Buser trouve  ses exemples de surfaces poss\'edant
 ``beaucoup'' de valeurs propres $\displaystyle \leq \frac 14$.
Soit $S$ une surface hyperbolique
de genre $g$ munie d'une d\'ecomposition en $2g-2$ pantalons $P_i$ \`a bord g\'eod\'esique dont toutes les courbes du bord ont une longueur 
inf\'erieure \`a $l$. Buser montre alors que les $2g-3$ premi\`eres valeurs propres de $S$ sont plus petites que $\epsilon (l)$ o\`u $\epsilon (l)$ tend vers $0$ avec $l$. En effet, quand $l$ est 
suffisamment petit, la premi\`ere valeur propre pour le probl\`eme de Dirichlet sur chacun de ces pantalons est inf\'erieure \`a $ \epsilon (l)$, o\`u $\epsilon (l)\to 0$ avec $l$. D\'efinissons pour $i=1,\cdots, 2g-2$,
 une fonction $f_i$ sur $S$   de la mani\`ere suivante :  sa restriction \`a   $P_i$ est  la premi\`ere fonction propre du probl\`eme de Dirichlet et $f_i$  est identiquement nulle   sur  le 
 compl\'ementaire de $P_i$.
On v\'erifie que   $f_i$
est dans le domaine de d\'efinition du Laplacien sur $S$\;;  son
quotient de Rayleigh, 
$\displaystyle \frac {\int _S\Vert \nabla f_i\Vert ^2 }{\int _S f_i^2}$  
\'egal \`a  la valeur propre de $f_i\vert S$ est inf\'erieur \`a $\epsilon (l)$.
Les fonctions $f_i$ forment un base ortho\-gonale
d'un espace vectoriel de dimension $2g-2$.
  La caract\'erisation variationnelle  (par mini-max)  des  valeurs propres,
  entra\^{\i}ne alors 
   $\displaystyle \lambda _{2g-3}\leq \epsilon (l)$.  
   Cette construction produit donc des exemples de surfaces pour lesquelles les $(2g-3)$ premi\`eres valeurs propres   sont petites.  
   D'autre part, un   r\'esultat de  Burton Randol,  qui \'etait \'egalement connu de Buser,  dit  qu'il existe une constante $c (g)>0$   telle que si dans une surface hyperbolique de genre $g$,   les      $(2g-3)$-premi\`eres valeurs propres sont inf\'erieures  \`a $c (g)$, alors 
  $\displaystyle \lambda _{2g-2}>\frac 14$ \cite {R2}. 
 Ce dernier r\'esultat a \'et\'e g\'en\'eralis\'e par R. Schoen, S. Wolpert et S. T. Yau : la valeur propre
 $\lambda _{2g-2}$ est comprise entre deux constantes positives qui ne d\'ependent que des longueurs des g\'eod\'esiques ferm\'ees \cite{SWP}.  
Signalons que dans le cas o\`u $g=2$, l'in\'egali\-t\'e   $\displaystyle \lambda _2 \geq \frac 14$
 a \'et\'e ob\-tenue par Paul Schmutz \cite {Sc}; dans cet article, la conjecture $\displaystyle \lambda _{2g-2}\geq \frac 14$ est aussi formul\'ee.  

\

  Nous \'etudierons ensuite le cas d'une surface hyperbolique d'aire finie.
   Du point de vue de la topologie,   une telle surface  $S$ est {\it 
de type $(g,n)$}, c'est-\`a-dire qu'elle est hom\'eomorphe au compl\'ementaire  de $n$ points dans une surface compacte  $\bar S$
de genre $g$  ; ces points, les {\it pointes}, ont un voisinage isom\'etrique 
au quotient du demi-plan hyperbolique $\Bbb H=\{\Im z >0\}$  par la translation $z\mapsto z+1$, la projection d'une droite horizontale sur $S$ est appel\'ee  un {\it horocycle}.
 Le spectre (de l'extension de Friedrichs) de l'op\'erateur $-\Delta$   sur $S$ est form\'e d'une partie discr\`ete et d'une partie continue contenue dans l'intervalle $\displaystyle [\frac 14,\infty[$ (provenant 
 des s\'eries d'Eisenstein) \cite{I}. 
 Le spectre discret 
$0=\lambda _0<\lambda _1\leq \cdots \leq \lambda _j\leq \cdots $
 est la r\'eunion (avec multiplicit\'e)
 du {\it   spectre parabolique},  c'est-\`a-dire de l'ensemble des valeurs propres $^p\lambda _j$
dont les fonctions propres   ont leurs moyennes nulles sur tous les horocycles centr\'ees aux pointes
et du {\it   spectre r\'esiduel}, un  ensemble  fini
de  valeurs propres $\displaystyle ^r\lambda _j \leq \frac 14 $ dont les fonctions propres  
apparaissent comme
r\'esidus du prolongement analytique des s\'eries d'Eisenstein. 
Le spectre parabolique n'est  pas  bien compris en g\'en\'eral\;; 
instable comme fonction de la m\'etrique, il pourrait bien m\^eme \^etre vide g\'en\'eriquement. 
Mais les valeurs propres inf\'erieures \`a $\displaystyle \frac 14$ sont en nombre fini \cite [p.112]{I}.
 Il est connu      que lorsque le genre de la surface compl\'et\'ee $\bar S$ vaut $0$ ou $1$, il n'existe pas de valeurs propres paraboliques   $\displaystyle \leq  \frac 14$ \cite {Hu}, \cite [Prop. 2] {O}.

 D'apr\`es
  \cite [Prop. 3]{O},  pour tout $\displaystyle \lambda \leq \frac 14$,   la multiplicit\'e  de $\lambda$ dans le spectre parabolique
est au plus $2g-3$; cette Proposition donnait  aussi que la multiplicit\'e de $\lambda$ (dans la r\'eunion des spectres 
r\'esiduels  et paraboliques)
est au plus $2g-3+n$, par la m\^eme m\'ethode que dans le cas compact.
Ici, nous allons montrer la g\'en\'eralisation suivante de ce r\'esultat.
\begin{theo}\label{2g-2parabolique}Soit $S$ une surface  hyperbolique
d'aire finie et de type $(g,n)$.
Alors
la $(2g-2+n )$-i\`eme valeur propre  
$\lambda _{2g-2+n }$    
est $\displaystyle >\frac 14$.
\end{theo}  
Il est possible  que le spectre parabolique v\'erifie la conclusion  du Th\'eor\`eme  \ref{2g-2}, c'est-\`a-dire que l'on ait  $\displaystyle ^p\lambda _{2g-2}>\frac 14$, mais nous n'avons pas pu l'\'etablir avec nos m\'ethodes topologiques\;; nous discutons \`a  la fin de l'article quelques questions d'analyse connexes.
 
 \
 
 \noindent{\bf{M\'ethode.}}
 Dans \cite {Bu2}, \cite {R2}, \cite{Sc}, \cite {SWP},  les minorations de valeurs propres  sont bas\'ees sur des arguments variationnels \cite{Bu2, Cha}.
 La    d\'emonstration du Th\'eor\`eme \ref{2g-2}
 repose sur une m\'ethode nouvelle  introduite par Bruno S\'evennec  dans
   \cite {Se}. Dans  cet  article, S\'evennec \'etudie la multiplicit\'e de la deuxi\`eme valeur propre des op\'erateurs de Schr\"odinger     sur les surfaces compactes et obtient le r\'esultat remarquable  suivant. 
\begin{theo}\cite[Theo. 4] {Se} Soit $S$ une surface riemannienne compacte de genre $g$.  Alors, pour toute fonction   lisse   $v :S\to \RR$, la multiplicit\'e de la deuxi\`eme 
  valeur propre   de $-\Delta +v$ est inf\'erieure \`a $2g+3$.
 \end{theo}
Sa d\'emonstration
utilise le lemme topologique
  suivant. Soit $\mathcal E$ l'espace  des fonctions propres pour la 
 deuxi\`eme 
 valeur propre non nulle. Soit
 $\Bbb S (\mathcal E)$ la sph\`ere unit\'e de cet espace (pour un  norme quelconque) et soit $\Bbb P (\mathcal E)$ l'espace projectif  sur $\mathcal E$, vu comme le quotient de $\Bbb S (\mathcal E)$  par l'involution 
 antipodale $\tau : f\mapsto \tau (f)=-f$ .   Supposons que  $\Bbb S (\mathcal E)$  admette une partition en $k$ ensembles 
 $\mathcal S_i$ tels que 
 \begin{enumerate}
 \item pour tout $i$, $\mathcal S_i$ est invariant par l'involution
  antipodale $\tau $ et 
 \item le rev\^etement $\mathcal S_i \to \mathcal S_i/\tau$ est trivial.
 \end{enumerate}
Alors, 
la dimension de 
$\Bbb S (\mathcal E)$ 
  est inf\'erieure \`a $k-1$, autrement dit 
     $\text{dim}(\mathcal E)\leq k$\;: 
c'est un r\'esultat  de type Borsuk-Ulam \cite [Lemme 8] {Se}.
 
 S\'evennec d\'efinit une partition  de $\Bbb S(\mathcal E)$  de la mani\`ere suivante. Pour toute fonction   $f\in \Bbb S(\mathcal E)$, les ouverts $S^+(f) =\{f>0\}$ et $S^-(f)=\{f<0\}$  ne sont pas vides.  Leurs premiers nombres de Betti v\'erifient $ b^1(S^+(f))+b^1(S^-(f))\leq b^1(S)$ \cite [Lemma 12] {Se}.  Les  atomes 
 $\mathcal S_i$
de la partition sont alors d\'efinis par  
$b^1(S^+(f))+b^1(S^-(f)) =i$\;; par d\'efinition, chaque atome
est   invariant par $\tau$ et le nombre d'atomes est au plus  \`a $b^1(S)+1$. S\'evennec montre que le rev\^etement $\mathcal S_i\to \mathcal S_i/\tau$ est trivial pour $i>0$. Le rev\^etement pour $i=0$  demande une \'etude particuli\`ere mais le lemme de type Borsuk-Ulam permet toutefois de conclure.

Dans \cite [Proposition 3] {O},   la  m\^eme m\'ethode 
est utilis\'ee  pour  minorer la multiplicit\'e d'une valeur propre $\lambda$  sur une surface $S$ de courbure n\'egative et   lorsque $\lambda$ est  inf\'erieur au bas du spectre $\lambda _0(\tilde S)$ du rev\^etement universel.  Quand  la valeur propre $\lambda $ n'est plus la premi\`ere valeur propre non nulle, les  ouverts 
$S^-(f)$
 et $S^+(f)$  ne sont pas n\'ecessairement connexes. Toutefois si $\lambda $ 
 est inf\'erieur \`a $\lambda _0(\tilde S)$,  ces   ouverts sont {\it incompressibles} (le  groupe fondamental de chacune de leur composante connexe s'injecte dans celui de $S$) et leur caract\'eristique d'Euler est strictement n\'egative \cite [Lemma 1] {O}.  L'auteur d\'efinit  alors une  partition de $\Bbb S (\mathcal E)$  en  les ensembles $\mathcal S_i$  
 d\'efinis par $\chi (S^+(f))+\chi (S^-(f)) =i$.  \`A cause de la propri\'et\'e 
d'incompressibilit\'e 
pr\'ec\'edente, l'entier $i$ peut prendre des valeurs entre $\chi (S)=2-2g$ et $-2$\;; donc la partition de $\Bbb S (\mathcal E)$ est form\'ee de $2g-3$ atomes au plus.  Chaque atome 
 $\mathcal S_i  $ est   invariant par $\tau$ et
  le rev\^etement $\mathcal S_i \to \mathcal S_i/\tau$ est trivial (pour ce dernier point, l'incompressibilit\'e pr\'ec\'edemment \'etablie permet de simplifier l'argument correspondant  dans \cite {Se}).
  Il obtient alors :  
 \begin{prop}\cite [Proposition 3]{O} Soit $S$ une surface riemannienne 
 de courbure n\'egative et soit $\lambda >0$ une  valeur propre du Laplacien
  $\leq \lambda _0 (\tilde S)$, le bas du spectre du rev\^etement universel de $S$. Alors la multiplicit\'e de $\lambda$ est inf\'erieure \`a $2g-3$.
  \end{prop}

Dans l'article pr\'esent, nous appliquerons toujours le m\^eme lemme pour montrer le th\'eor\`eme  \ref{2g-2}.  
  Pour cela, notons 
 $\mathcal E$
  l'espace  vectoriel 
   engendr\'e par les fonctions $\lambda$-propres pour $\displaystyle \lambda\leq \frac 14$ 
   (fonctions constantes comprises). D\'emontrer
le  th\'eor\`eme  \ref{2g-2} revient  \`a montrer que la dimension $m$ de $\mathcal E$ est inf\'erieure \`a $2g-2$ (remarquons qu'un cas particulier 
du th\'eor\`eme  \ref{2g-2} est
 la proposition 4   rappel\'ee 
  ci-dessus
puisque $\mathcal E$ contient l'espace des fonctions constantes). Les articles \cite {Se} et \cite {O} reposaient en partie sur la connaissance de  la structure de {\it l'ensemble nodal} -- c'est-\`a-dire l'ensemble des z\'eros -- d'une fonction propre\;:
d'apr\`es un r\'esultat classique de Cheng, cet ensemble est localement diff\'eomorphe par un diff\'eomorphisme ambiant au lieu des z\'eros d'un poly\-n\^ome harmonique \cite {Che}.
Pour une combinaison  lin\'eaire  de fonctions propres, la si\-tuation est 
  plus compliqu\'ee.
 C'est pour cette raison  que l'une des   hypoth\`eses du th\'eor\`eme \ref{2g-2}
 est que  la m\'etrique est analytique. Les fonctions propres sont  alors analytiques\;; une  combinaison lin\'eaire finie $f=\sum f_j$ l'est
 aussi et on peut \'etudier son ensemble nodal comme celui d'une fonction analytique. Un argument de Lojasiewicz montre alors 
que l'ensemble nodal de $f$
est la r\'eunion  d'un ensemble discret
et d'un graphe plong\'e dans $S$, le {\it graphe nodal}.
 La premi\`ere section de l'article est consacr\'ee \`a l'\'etude de 
{\it l'ensemble nodal} des fonctions de $\mathcal E$\;: on red\'emontre en particulier ce r\'esultat de Lojasiewicz. 

Dans la section suivante, on d\'efinit {\it les surfaces  caract\'eristiques de $f$}, not\'ees
$\Sigma ^+ (f)$
 et $\Sigma ^- (f)$. Ces surfaces joueront le m\^eme r\^ole que les surfaces $S^+ (f)$ et $S^- (f)$  ci-dessus.
 On remarque en effet que les surfaces $S^\pm (f)$, telles qu'elles \'etaient d\'efinies auparavant  ne  sont pas toujours 
    incompressibles  et perdent d\`es lors la propri\'et\'e  de stabilit\'e qui permettaient de v\'erifier  que les rev\^etements $\mathcal S_i\to \mathcal S_i /\tau$ \'etaient triviaux\;; on  les remplace par des surfaces  $\Sigma ^\pm (f)$ construites de la mani\`ere suivante.
  On  \'elimine d'abord du graphe nodal de $f$ les composantes connexes qui sont   contenues dans des disques. Les  composantes du 
  compl\'ementaire du graphe obtenu sont incompressibles par construction. De plus,  sur chaque composante connexe  $f$   garde un signe constant  sauf peut-\^etre sur une r\'eunion finie de disques (qui contiennent les composantes homotopes \`a $0$ que l'on a enlev\'ees).  On note 
  $S^\pm(f)$
  la r\'eunion des composantes sur lesquels
 $f$ garde le m\^eme signe $\pm$ et qui ont une caract\'eristique d'Euler $<0$.  L'une de ces surfaces pourrait \^etre   vide\;; toutefois, en  consid\'erant  le  quotient de Rayleigh de $f$, on montre que l'une au moins des deux surfaces $S^\pm(f)$ a une caract\'eristique d'Euler   strictement n\'egative. On d\'efinit alors  $\Sigma^\pm (f)$ comme la r\'eunion d'un c\oe ur compact  de $S^\pm (f)$ et des anneaux \'eventuels, qui joignent deux composantes de $S^\pm(f)$.

 Dans la troisi\`eme section, nous d\'emontrons le th\'eor\`eme \label{2g-2} en adaptant la d\'emonstration de \cite {O}. On d\'efinit  une partition de $f$ en ensembles $\mathcal S_i$  caract\'eris\'es par la propri\'et\'e que $\chi  (\Sigma ^+(f))
 +     \chi  (\Sigma ^-(f)) = i$.
 Cette fois, les valeurs prises par $i$ vont de $2-2g $ \`a $-1$; il y a donc au plus $2-2g$  valeurs. Par d\'efinition, chaque atome est invariant par  l'involution $\tau$.  Le lemme de  type Borsuk-Ulam permet de conclure. 
  
Dans la section suivante, nous  adaptons l'argument au cas des surfaces d'aire finie. 
 
Nous remercions le referee dont les remarques ont permis  
 d'am\'eliorer la pr\'esentation de l'article.
  
\
 
\section{L'ensemble nodal d'une somme finie de  fonctions propres}
 
Soit   $f=\sum  _1^n f_j$ avec $\Delta f_j+\lambda _j f_j=0$ une combinaison finie de fonctions $\lambda_j$-propres.
L'ensemble des z\'eros de $f$, $f^{-1}\{0\}$ 
est appel\'e {\it l'ensemble nodal de $f$}; on le note 
$Z (f)$.
Nous nous int\'eressons ici \`a la topologie de 
$Z(f)$. 
Lorsque $S$ est 
munie d'une m\'etrique $C^\infty$,
un th\'eor\`eme de Cheng \cite{Che} d\'ecrit la topologie de l'ensemble nodal   d'une fonction propre\;: pour tout $p\in Z(f)$, l'intersection de $Z(f)$ avec un petit voisinage de $p$ est hom\'eomorphe \`a l'ensemble nodal
d'une fonction harmonique, c'est-\`a-dire \`a l'ensemble des z\'eros de   $\displaystyle \Re (z^k)$.  Pour une m\'etrique de courbure $-1$, plus g\'en\'eralement pour une m\'etrique analytique, nous allons montrer que la topologie de l'ensemble nodal  d'une combinaison lin\'eaire $f$ n'est gu\`ere plus compliqu\'e.
 \begin{prop}\label{lieunodal} Soit $S$ une surface  
  munie d'une m\'etrique Riemannienne analytique. 
Soit $f:S\to \RR$ une combinaison lin\'eaire non nulle d'un nombre fini de fonctions propres du Laplacien sur $S$.
Alors l'ensemble nodal $Z(f)$   est la r\'eunion localement finie d'un graphe sans sommets libres et de points isol\'es.
\end{prop}
\noindent{\bf{D\'emonstration.}} Soit $p\in Z(f)$. Si le vecteur gradient $\nabla f (p)$ n'est pas nul, alors $Z(f)$ intersecte un voisinage de $p$ comme une sous-vari\'et\'e de dimension $1$ d'apr\`es le th\'eor\`eme des fonctions implicites\;; ceci 
d\'emontre la proposition  dans ce cas.

 Puisque la m\'etrique est analytique, les fonctions propres du Laplacien sont solutions d'une \'equation elliptique \`a coefficients analytiques\;;  ce sont donc des fonctions analytiques.  Par cons\'equent,
   $f$ est analytique. Un  th\'eor\`eme fondamental de Lojasiewicz \cite{Lo}  donne  une stratification du  lieu des z\'eros d'une fonction analytique $\RR^n\to \RR$. En dimension $2$, l'argument est plus simple et nous le rappelons bri\`evement. 
 Soit $p\in Z (f)$. Le Th\'eor\`eme de pr\'eparation de Weierstrass permet d'\'ecrire $f$
sur un voisinage $V(p)$ de $p$, param\'etr\'e par le carr\'e $]-\delta, \delta [\times ]-\delta ,\delta[$,
comme un produit 
$f(x,y)=Q(x,y)H(x,y)$\;:
ici $Q$ 
est  une fonction analytique
qui ne s'annule pas  et 
$H$ est  un  {\it polyn\^ome adapt\'e}, c'est-\`a-dire un polyn\^ome
de la forme
 $$
H(x,y)= y^m+a_1 (x)y^{m-1}+\cdots +a_m (x)
$$
pour des fonctions  $a_i(x)$
analytiques sur $]-\delta, \delta[$
et  qui v\'erifient $a_i(0)=0$. 
Notons
$\mathcal A$ l'anneau   des fonctions \`a valeurs r\'eelles analytiques en la variable $x\in ]-\delta, \delta  [$. 
Ainsi $H$ est un polyn\^ome dans l'anneau $\mathcal A [y]$ des polyn\^omes en la variable $y$.

Une cons\'equence de cette \'ecriture est que l'ensemble des points de $]-\delta, \delta [\times ]-\delta ,\delta[$ o\`u   $f$  s'annule
est exactement celui   o\`u $H$ s'annule. L'anneau $\mathcal A$ est factoriel, et donc   $\mathcal A[y]$ l'est aussi. En raisonnant comme dans \cite{Lo}, on se ram\`ene au cas o\`u le discriminant de $H$ est non nul tout 
 en conservant la propri\'et\'e que le lieu des z\'eros de $f$ est aussi celui de $H$. Pour cela on d\'ecompose le polyn\^ome $H$ en produit de facteurs 
 irr\'eductibles 
   $H=H_1^{r_1}...H_k^{r_k}$. Le nouveau polyn\^ome $H'=H_1...H_k$ a exactement les m\^emes z\'eros  que $H$ 
    dans $V$  et son discriminant, un \'el\'ement de l'anneau $\mathcal A$ est non nul.
 Si   ce  discriminant ne s'annule pas en $0$,  alors  $0$ est racine simple du polyn\^ome $H(0,y)$ et on peut suivre cette racine  comme  une   fonction  analytique 
 $\zeta (x)$. Supposons que ce discriminant s'annule en $0$ ; alors puisque le discriminant est une fonction analytique de $x$, sa valeur au point
  $x$   est   non nulle
  pour tout $x\neq 0$ $\vert x \vert <\epsilon$ si $\epsilon$ est suffisamment petit. Sur chacun des   intervalles $]0,\epsilon[$ et $]-\epsilon, 0[$, on peut suivre les racines de $H(x,.)$ comme les graphes de fonctions analytiques $\zeta _1 (x),\cdots, \zeta _m(x)$.  Parmi ces racines, il y a celles qui sont r\'eelles (et distinctes)
$\zeta _1, \cdots,\zeta _k$; les autres sont imaginaires conjugu\'ees deux-\`a-deux. Lorsque $x\to 0^\pm$, la limite des fonctions $\zeta _i (x)$ est   n\'ecessairement $0$. On voit ainsi $Z(f)=Z(H')$ localement comme la r\'eunion des graphes (au-dessus de $]-\epsilon, 0]$ et de $[0,\epsilon[$)
des fonctions $\zeta _i$  (si   ces racines ne sont pas toutes imaginaires) ou comme r\'eduit au point $p$   (lorsque toutes les racines sont complexes conjugu\'ees deux-\`a-deux). On v\'erifie que le nombre de ``branches'' en $p$ 
est pair (c'est un cas tr\`es particulier d'un th\'eor\`eme de Sullivan).
Ceci montre  la Proposition.
 \hfill\qed\enddemo
Quand $S$ est compacte, 
l'ensemble nodal est la r\'eunion d'un graphe compact et d'un ensemble fini de z\'eros isol\'es.
 Quand $S$  n'est pas compacte, on ne dispose que de  la description locale ; a priori il pourrait y avoir un nombre infini de z\'eros   contenus
dans une pointe de $S$. Dans le cas d'une surface hyperbolique,
   l'ensemble nodal 
d'une fonction $\lambda$-propre   avec $\displaystyle \lambda \leq \frac 14$,
 a  une adh\'erence dans la 
   surface compl\'et\'ee $\bar S$ qui  est un graphe compact
   \cite [Lemma 6]{O}. Nous ne savons pas si 
     un r\'esultat analogue est v\'erifi\'e
par  l'ensemble nodal
 d'une fonction dans $\mathcal E$.
   
 \
 
\section
{La surface caract\'eristique d'une somme finie  de fonctions propres}

Dans cette section, $S$ est  une surface compacte munie d'une m\'etrique analytique de courbure n\'egative.
On note  $\mathcal E$ le $\RR$-espace vectoriel
  engendr\'e par les fonctions $\lambda$-propres, pour $\lambda\leq  \lambda _0(\tilde S)$.
 Dans cette section, nous associons \`a toute fonction non nulle
 $f\in \mathcal E$ une sous-surface de $S$, compacte et incompressible,
 {\it la surface  caract\'eristique de $f$}.

On d\'efinit {\it le graphe nodal}   d'une fonction $f\in \mathcal E$ comme la r\'eunion, not\'ee
$\gamma (f)$, des composantes connexes de son ensemble nodal qui sont des graphes. Le graphe $\gamma (f)$ peut
\^etre vide. 
Toutefois il ne l'est pas  lorsque $f$ change de signe.
Sur chaque composante connexe de $S-\gamma (f)$, on remarque que  le signe de $f$ est bien 
d\'efini puisque ses z\'eros y forment un ensemble discret.
 On note $G(f)$ le graphe obtenu en retirant de $\gamma (f)$
les composantes connexes qui sont contenues dans des disques. 
Chaque composante connexe de $S-G(f)$  est \'egale \`a la r\'eunion d'une composante de $S-\gamma (f)$ et d'un nombre fini de disques contenus dans $S$ et disjoints deux-\`a-deux. 
On peut donc la munir d'un signe $\pm$, \`a savoir celui  de la fonction $f$ sur  cette composante de 
$S-\gamma (f)$. 
On note $C^+(f)$ (resp. $C^- (f)$)
 la r\'eunion des composantes de $S-G(f)$ sur lesquelles le signe de $f$ est $+$ (resp. $-$).
Lorsque  $\gamma (f)$ est vide, la fonction $f$  garde le m\^eme
   signe\;; on d\'efinit alors $C^\pm(f)$  comme la surface $S$ avec  le signe \'evident.
 Dans tous les cas, $C^+ (f)$ ou $C^- (f)$
  est vide lorsque chaque composante de 
 $\gamma(f)$ est contenue dans  un disque.

Comme cons\'equence de cette construction,  les 
surfaces $C^+(f)$ ou de $C^-(f)$ sont
  {\it incom\-pres\-sibles}, c'est-\`a-dire que
le   groupe fondamental de chacune de leurs composantes connexes
  s'injecte dans $\pi _1 (S)$ par l'application 
 d'inclu\-sion.    La r\'eunion des composantes connexes de 
 $C^+(f)$ (resp. de $C^-(f)$) qui ne sont pas des disques ni des anneaux est une
 surface 
not\'ee  $S^+(f)$ (resp. $S^-(f)$)\;: elle est 
 peut-\^etre vide, peut-\^etre non connexe mais elle est incompressible.
  
 On note $\chi  ^+(f)$ 
 (resp. $\chi ^-(f)$ la caract\'eristique d'Euler de $S^+(f)$ (resp. de $S^-(f)$) (on convient ici que la caract\'eristique d'Euler d'une surface vide est $0$). L'incompressibilit\'e entra\^{\i}ne que $\chi ^+ (f)+\chi ^-(f)$ est sup\'erieur \`a $\chi (S)$. Par d\'efinition, on a $\displaystyle \chi ^\pm(f)\leq 0$ avec \'egalit\'e uniquement si $S^\pm(f)$ est vide.
\begin{cla} \label{moinsque1}Pour toute fonction non nulle $f\in \mathcal E$, on a : $\chi ^+ (f)+\chi ^-(f)<0$. 
\end{cla}
\noindent{\bf{D\'emonstration.}} Lorsque le graphe nodal $\gamma (f)$ est vide, $S^+ (f)=S$ ou $S^-(f)=S$
 selon que $f$ est $\geq 0$ ou $\leq 0$.  Il n'y a rien \`a montrer dans ce cas.

Dans le cas g\'en\'eral,  soit $\displaystyle f\in \mathcal E$
une fonction non nulle. 
Puisque  $f$ est   somme de fonctions propres pour des valeurs  propres $\leq \lambda _0 (\tilde S)$,  le quotient de Rayleigh $\displaystyle \frac {\int _S\Vert \nabla f \Vert ^2}
{\int _S f^2}$ est inf\'erieur ou \'egal \`a $\lambda _0 (\tilde S)$, avec \'egalit\'e  uniquement si $f$ est elle-m\^eme une fonction $\lambda _0 (\tilde S)$-propre. 
Notons $S_j$, pour $j=1,\cdots k$   les composantes de $C^+\cup C^-=S-G(f)$ ;
on a  $\displaystyle \int _S \Vert \nabla f \Vert ^2=\sum _j\int _{S_j}\Vert \nabla f \Vert ^2$ et 
$\displaystyle \int _S f^2=\sum _{j} \int _{S_j}f^2$. Donc pour   au moins un indice $j\leq k$, le quotient de Rayleigh  
$\displaystyle \frac {\int _{S_j}\Vert \nabla f \Vert ^2}{\int _{S_j} f^2}$ est inf\'erieur o\`u \'egal \`a $\lambda _0 (\tilde S)$. Montrons  en raisonnant par l'absurde que la caract\'eristique d'Euler de   $S_j$ est strictement n\'egative. Dans le cas contraire, le rev\^etement $\tilde S _j$ de $S$ de groupe fondamental $\pi _1 (S_j)$ est ou bien le rev\^etement universel $\tilde S$, ou bien le quotient de ce dernier par un groupe cyclique. Il est bien connu que dans les deux cas, le bas du spectre du Laplacien est $\lambda _0 (\tilde S)$\;:  c'est la 
  d\'efinition m\^eme du bas du spectre lorsque ce rev\^etement est $\tilde S$,  \c ca d\'ecoule d'un th\'eor\`eme de Brooks  \cite {Br} lorsque c'est un quotient cyclique (dans le cas des surfaces qui nous int\'eresse, on peut le voir aussi d'une mani\`ere plus   \'el\'ementaire  en consid\'erant les quotients de Rayleigh). La surface $S_j$ se rel\`eve dans $\tilde S _j $
et la fonction $f\vert S_j$ prolong\'ee par $0$ dans le compl\'ementaire de ce relev\'e est dans le domaine du Laplacien sur $\tilde S _j$\;: bien que $f\vert S_j$
 n'est pas d\'erivable,  sa d\'eriv\'ee au sens des distributions est dans $L^2$. Le quotient de Rayleigh de $f\vert S_j$ est   $\leq  \lambda _0 (\tilde S)$\;: on a donc une contradiction (dans le cas o\`u ce quotient est \'egal \`a $\lambda _0 (\tilde S)$, on utilise  le fait que si le quotient de Rayleigh d'une fonction est \'egal au bas du spectre du Laplacien, alors cette fonction est une fonction
 propre).
\hfill\qed\enddemo 
 
\

\noindent{\bf{D\'efinition.}}  Pour chaque composante   $S_j ^\pm(f)$
de $S^\pm(f)$ de caract\'eristique d'Euler n\'egative, choisissons un {\it c\oe ur compact}, c'est-\`a-dire une surface compacte 
$K_j^\pm (f)\subset S_j ^\pm(f)$ telle que l'inclusion soit une \'equivalence d'homotopie.  On d\'efinit  
 la surface    $\Sigma ^+(f)$   (resp. $\Sigma ^-(f)$)  comme \'etant
 la r\'eunion des c\oe urs compacts $K_j^+(f)$ (resp. $K_j ^-(f)$) et des  composantes  \'eventuelles  du compl\'ementaire $S-\bigcup K_j^+(f)$, (resp. $S-\bigcup K_j ^-(f)$)  qui sont des anneaux. 
 Donc $\Sigma ^+ (f)$ (resp. $\Sigma ^- (f)$) s'obtient \`a partir de $\bigcup K_j ^+$ (resp.
 $\bigcup K_j ^-$ ) en lui rajoutant, s'il en existe, les anneaux   entre les composantes du bord de  $\bigcup K_j ^+(f)$  (resp. de $\bigcup K_j ^-(f)$) qui sont homotopes.
 La surface   $\Sigma (f)=\Sigma ^+ (f)\cup \Sigma ^-(f)$ est   appel\'ee   {\it la surface caract\'eristique de $f$}, et
 $\Sigma ^+ (f)$ (resp.    $\Sigma ^-(f)$) {\it la surface  caract\'eristique    positive} (resp. {\it n\'egative}). La d\'efinition de
 ces surfaces d\'epend uniquement du choix des c\oe urs compacts : elles sont donc bien 
 d\'efinies \`a isotopie pr\`es.
Par construction, la  caract\'eristique d'Euler  de $\Sigma ^+(f)$ vaut $\chi ^+ (f)$, celle 
de $\Sigma ^-(f)$ vaut $\chi ^-(f)$. 
On a bien entendu $\Sigma ^+(-f)=\Sigma ^- (f)$ et $\Sigma ^-(-f)=\Sigma ^+ (f)$.

\
 
\section{D\'emonstration du th\'eor\`eme \ref{2g-2}}

Soit $m$ la dimension de l'espace vectoriel $\mathcal E$ introduit dans la section pr\'ec\'edente.
Le th\'eor\`eme \ref{2g-2} d\'ecoulera de l'in\'egalit\'e  $m\leq 2g-2$.
Notons $\Bbb S(\mathcal E)$  la sph\`ere unit\'e de   
$\mathcal E$  (pour une norme arbitraire) 
et $\Bbb P(\mathcal E)$  l'espace projectif sur $\mathcal E$, quotient de $\Bbb S(\mathcal E)$ par l'involution    $f\mapsto -f$.

 Pour chaque entier $i$ avec $2-2g\leq i\leq -1$, on note
$$
\mathcal S_i=\{f\in \Bbb S (\mathcal E)\quad\vert \quad  \chi ^+ (f)+\chi ^- (f)=i\}.
$$
 D'apr\`es l'affirmation \ref{moinsque1},  $\Bbb S  (\mathcal E)=\bigcup _{2-2g}^{-1}\mathcal S_i$.
D'autre part, chaque ensemble $\mathcal S_i$ est inva\-riant par l'involution antipodale.
L'espace projectif $\PP (\mathcal E)$  est r\'eunion des  ensembles $\mathcal P_i$, quotients de $\mathcal S_i$ par l'involution antipodale.
\begin{lem}Pour tout entier $i$, $2-2g\leq i\leq -1$, le rev\^etement $\mathcal S_i\to \mathcal P_i$ est trivial.
\end{lem}
 \noindent{\bf{D\'emonstration.}}
Soit $f\in \mathcal S_i$.  Reprenons les notations introduites dans  la d\'efinition  de la surface 
 {caract\'eristique} de $f$ : $K_j^\pm(f)$ est un c\oe ur compact de $S_j^\pm(f)$, composante connexe de caract\'eristique d'Euler n\'egative de $S^\pm (f)$. On supposera ce c\oe ur compact choisi de sorte que les composantes connexes de $\gamma (f)$ qui seraient \'eventuellement contenues dans $S_j ^\pm (f)$ soient contenues dans l'int\'erieur de $K_j ^\pm (f)$.

Pour toute fonction $g\in \mathcal E$, suffisamment proche de $f$,  $K_j^+(f)$ est alors contenu dans une composante connexe  $S_l^+(g)$
de $S^+ (g)$. 
Fixons un voisinage $V(f)$  de $f$ dans $\mathcal S_j^\pm(f)$  tel que  ces inclusions aient lieu pour  chaque  surface $K_j^\pm (f)$. 
Nous allons montrer que pour toute fonction $g\in \mathcal S_i \cap V(f)$, les surfaces caract\'eristiques  $\Sigma ^+(f)$ et   $\Sigma ^+(g)$  (resp. $\Sigma ^-(f)$ et   $\Sigma ^-(g)$) sont isotopes.
Choisissons maintenant les c\oe urs compacts $K_l^\pm(g)$ 
des surfaces $S^\pm (g)$ de sorte que lorsque $K_j^\pm (f)$ est contenu dans $S_l ^\pm (g)$, il soit
 contenu dans l'int\'erieur de 
$K_l^\pm (g)$. 
Supposons que deux composantes du bord des surfaces 
$K_j^+ (f)$ sont homotopes dans $S$\;; alors l'anneau qui r\'ealise l'homotopie est contenu dans $\Sigma ^+ (f)$, par d\'efinition de la surface caract\'eristique\;;  puisque  cet anneau joint deux courbes de $K^+ (g)$, il est contenu dans l'une des composantes connexes de $\Sigma ^+ (g)$ toujours par d\'efinition de la surface caract\'eristique. On en d\'eduit que chaque composante connexe de $\Sigma ^\pm (f)$ est contenue dans une composante connexe de 
 $\Sigma ^\pm (g)$ (de m\^eme signe). 
 Puisque $\Sigma ^+(f)$ et $\Sigma ^- (f)$ sont incompres\-sibles dans $S$, elles sont incompres\-sibles dans $\Sigma ^+(g)$ et $\Sigma ^- (g)$ respectivement. 
 En particulier, leurs caract\'eristiques d'Euler v\'erifient  $\chi ^+(f)\geq \chi ^+ (g)$  et  $\chi ^- (f)\geq \chi ^- (g)$ ;
de plus ces in\'egalit\'es sont des  \'egalit\'es si et seulement si les surfaces $\Sigma ^+ (f)$ et   $\Sigma ^+ (g)$ sont isotopes ainsi que  $\Sigma ^- (f)$ et $\Sigma ^- (g)$. 
Or, puisque $g\in \mathcal S_i$, on a 
l'\'egalit\'e  $\chi ^+(f)+\chi ^-(f)
 =\chi ^+(g)+\chi ^-(g)$.
 Donc $\Sigma ^+(f)$ et $\Sigma ^+(g)$
sont isotopes, de m\^eme que 
$\Sigma ^- (f)$ et $\Sigma ^-(g)$. 
 
 Puisque les   fonctions ``classe d'isotopie de $\Sigma ^+(f)$'' et  
 ``classe d'isotopie de  $\Sigma ^-(f)$''  sont localement constantes sur $\mathcal S_i$, elles
 sont constantes sur les composantes connexes. Les fonctions $f$ et $-f$ ne peuvent pas \^etre dans la m\^eme composante connexe de $\mathcal S_i$, car  sinon  $\Sigma ^+(f)$ et $\Sigma ^-(f)$ seraient isotopes
 d'apr\`es le paragraphe pr\'ec\'edent. Or deux surfaces de caract\'eristique d'Euler n\'egative contenues dans $S$,  disjointes et incompressibles ne peuvent pas \^etre isotopes. Le rev\^etement $\mathcal S_i\to \mathcal P_i$ est donc le
  rev\^etement trivial.
\hfill\qed\enddemo 
 
  Le rev\^etement  double 
  $\Bbb S (\mathcal E)\to \Bbb P(\mathcal E)$  est
  d\'ecrit par une classe $\alpha \in H^1(\Bbb P(\mathcal E),\ZZ/2\ZZ)$. Chaque rev\^etement $\mathcal S_i\to \mathcal P_i$ est d\'ecrit par la classe de cohomologie de \v{C}ech, $\alpha \vert _{\mathcal P_i}$. Puisque  chacun de ces rev\^etements est le rev\^etement trivial, on a $\alpha \vert _{\mathcal P_i}=0.$ Puisque $\Bbb P(\mathcal E)$ est la 
  r\'eunion des  ensembles $\mathcal P_i$ et que ceux-ci sont au plus au nombre de $2g-2$, on a : $\alpha ^{2g-2}=0$
  (cf.  \cite  
  [Lemme 8] {Se}). Puisque
   $\alpha$ est d'ordre $m$ dans l'anneau de $\ZZ/2\ZZ$-cohomologie
   de $\Bbb P(\mathcal E)$, on a 
  $ m\leq 2g-2$.
 \hfill\qed\enddemo
  \
  
\section{D\'emonstration du Th\'eor\`eme \ref{2g-2parabolique}}

Nous ne ferons qu'indiquer les modifications \`a apporter au raisonnement pr\'ec\'edent.  Ici $S$ est une surface hyperbolique d'aire finie.
On note $\bar S$ la surface compacte  de genre $g$ obtenue en rajoutant \`a $S$ ses pointes. 
On note   $m$ la dimension de l'espace vectoriel $\mathcal E$
engendr\'e par les fonctions $\lambda$-propres   v\'erifiant $\lambda \leq \frac 14$. La sph\`ere unit\'e de $\mathcal E$, l'espace projectif sont not\'es $\Bbb S (\mathcal E)$,  $\Bbb P (\mathcal E)$.

Rappelons d'abord quelques propri\'et\'es   du comportement des fonctions propres au voisinage des pointes.
Pour $\displaystyle \lambda \leq \frac 14$, on sait que   l'ensemble nodal d'une fonction $\lambda$-propre intersecte 
 chaque bout comme une  r\'eunion d'arcs disjoints aboutissant \`a la pointe \cite [Proposition 6] {O}; \`a l'oppos\'e de ce comportement,
une fonction $\lambda$-propre qui n'est pas parabolique en une pointe, tend vers l'infini au voisinage de cette pointe (elle y est \'equivalente \`a $cy^s$ o\`u $s$ est d\'efini par  $\displaystyle s(1-s)=\lambda$, $\displaystyle s\geq \frac 12$)  et son ensemble nodal \'evite donc une voisinage de la pointe.
Cette derni\`ere propri\'et\'e 
 est  donc v\'erifi\'ee par les fonctions dans  $  \mathcal E$\;: si 
 $f\in \mathcal E$   n'est pas parabolique au voisinage d'une pointe, alors elle tend vers l'infini en cette pointe, puisque   le terme non parabolique,  en $y^s$, domine.  Par contre, si $f\in \mathcal E$ est somme de fonctions propres paraboliques, alors sa moyenne est  nulle sur chaque horocycle centr\'e en la pointe, 
    mais nous ne savons   pas comment en d\'eduire que la topologie de l'ensemble nodal $Z (f)$ est  aussi simple pr\`es de la pointe que celle de l'ensemble nodal d'une fonction propre. On peut uniquement dire que l'ensemble nodal s'accumule sur la pointe (puisque sur chaque horocycle centr\'e  en la pointe, $f$ \'etant de moyenne nulle,  s'y 
  annule en deux points au moins). 
Toutefois, les fonctions  dans $\mathcal E$  
sont analytiques et leur ensemble nodal est donc d\'ecrit par la Proposition   \ref{lieunodal}\;: c'est la r\'eunion d'un graphe nodal $\gamma (f)$ et d'un ensemble discret de $S$.  

Les surfaces $S^\pm (f)$ et $\Sigma ^\pm(f)$
 sont d\'efinies  de la mani\`ere suivante. On note $G(f)$ le graphe 
obtenu \`a partir de ${\gamma(f)}$ en retirant les composantes connexes 
born\'ees qui sont  homotopes \`a $0$  dans  $ S$.
En utilisant que  les points d'accumulation
de $\gamma(f)$ dans $\bar S$ sont  contenus dans l'ensemble fini des pointes, on  construit un ferm\'e de $S$,  r\'eunion de disques ferm\'es deux-\`a-deux disjoints, contenus dans les composantes connexes de $  S-G(f)$
et qui contiennent la r\'eunion des   composantes born\'es de $\gamma (f)$ qui sont homotopes \`a 
 $0$ dans $S$. 
On munit chaque composante connexe de $ S- {G(f)}$ d'un signe, 
d\'efini  de la m\^eme mani\`ere que  dans le cas compact.
Contrairement au cas compact, on ne sait pas si $S-G(f)$ a un nombre fini de composantes connexes.
On d\'efinit 
alors $S^+(f)$ (resp. $S^-(f)$)
comme  la r\'eunion des composantes connexes de 
$  S-  {G(f)}$  de caract\'eristique d'Euler n\'egative et de signe $+$ (resp.  de signe -);
et ensuite $\chi ^\pm  (f)$  comme la caract\'eristique d'Euler de $S^\pm (f)$.
Par construction $S^+(f)$ et $S^-(f)$ sont des surfaces incompressibles dans $S$.  Le graphe $G(f)$ est de mesure nulle;   le raisonnement utilis\'e dans l'affirmation \ref{moinsque1} donne 
alors que $\chi ^+(f)$ ou $\chi ^-(f)$ est inf\'erieur \`a $-1$.  
Puisque chacune de leurs composantes connexes    est incompressible et de caract\'eristique d'Euler n\'egative, 
 $S^+(f)$ et  $S^-(f)$ n'ont qu'un nombre fini de composantes connexes.
On construit alors la surface caract\'eristique $\Sigma ^\pm (f)$ comme dans le cas compact.
 La sph\`ere $\Bbb S(\mathcal E)$ est la r\'eunion 
 des ensembles $\mathcal S_i$ tels que $\chi ^+(f)+\chi ^-(f)=i$. Il y a,
 cette fois,
  au plus     $ \vert \chi (S)\vert = 2g-2+n$ ensembles $\mathcal S_i$ qui ne sont pas vides (puisque $S$ porte une m\'etrique hyperbolique, $2g-2+n >0$). Le raisonnement  utilis\'e dans le cas compact donne 
 $m \leq 2g-2+n$\;: c'est le r\'esultat cherch\'e.
\hfill\qed\enddemo

\

\section{Conclusion}

La minoration du  Th\'eor\`eme \ref{2g-2parabolique} concerne la $(2g-2+n)$-i\`eme valeur propre du spectre discret,  r\'eunion des spectres r\'esiduel et parabolique.  Elle est optimale d'apr\`es le r\'esultat suivant.

\begin{remark} Pour tout couple d'entiers
$(g,n)$ avec $2g-3+n\geq0$, pour tout $\epsilon >0$,
il existe   une m\'etrique hyperbolique de volume fini 
et de type $(g,n)$ telle que $\lambda _{2g-3+n} \leq \epsilon$.
\end{remark}

\c Ca d\'ecoule de la construction de Buser
d\'ecrite dans l'introduction. Une   surface ho\-m\'eo\-morphe au compl\'ementaire de $n$ disques dans une surface de genre $g$
 poss\`ede une d\'ecomposition en
$2g-2+n$ pantalons $P_i$. 
On construit alors une m\'etrique hyperbolique sur $S$ en recollant le long du bord ces pantalons munis de m\'etrique  hyperbolique dont les composantes de bord correspondent \`a des pointes ou bien sont des g\'eod\'esiques de longueur $l$ (le nombre de pointes sur le pantalon $P_i$ d\'epend de la combinatoire de $P_i$ dans la d\'ecomposition). 
 Alors si $l$ est assez petit,    la surface hyperbolique obtenue a $2g-3+n$ valeurs propres inf\'erieures \`a $\displaystyle \frac 14$.

\
 
Mais dans cette construction,  il est difficile de savoir si les petites
 valeurs propres  sont  dans le spectre r\'esiduel ou dans le spectre parabolique. Un r\'esultat de Peter Zograf \cite {Z} garantit que si
  le nombre $n$ de pointes 
  est suffisamment grand par rapport \`a $g$, alors {\it toute} m\'etrique hyperbolique de volume fini et de type $(g,n)$  poss\`ede une valeur propre non nulle inf\'erieure \`a $\displaystyle \frac 14$.  Dans le m\^eme article, Zograf construit  aussi un exemple d'un sous-groupe d'indice fini (avec torsion) de $\text{PSL(2,}\ZZ)$ qui poss\`ede une valeur propre inf\'erieure \`a $\displaystyle \frac 14$ qui est dans le spectre parabolique. Dans ce cas, l'absence de valeurs propres  r\'esiduelles (autres que $0$ bien entendu) est garantie par leur absence dans le cas de $\text {PSL(2},\ZZ)$. 
Zograf montre ainsi que, pour cet exemple, la valeur propre inf\'erieure \`a $\displaystyle \frac 14$ qu'il a trouv\'ee   est n\'ecessairement   dans le spectre parabolique.
 
 \

Sur le spectre r\'esiduel,  le seul r\'esultat que nous connaissions et qui soit commun \`a toutes les surfaces de type $(g,n)$ est que la multiplicit\'e d'une valeur propre r\'esiduelle $\displaystyle \leq \frac 14$ est   $\displaystyle \leq n$\;; c'est une cons\'equence de l'existence du prolongement m\'eromorphe des s\'eries d'Eisenstein \cite [p. 112]{I}. 
  
Quant au spectre parabolique,    on sait que 
pour $\lambda \leq \frac 14$ la multiplicit\'e de $\lambda$ dans le spectre parabolique est major\'ee  
 par $\leq 2g-3$  \cite [ Prop. 2] {O},
 ind\'ependamment  donc du nombre de pointes  $n$. 
Nous pensons  qu'il existe une majoration analogue du  nombre total de valeurs propres parabo\-liques inf\'erieures \`a $\displaystyle \frac 14$. 
Notons $0<^p\lambda _1 \leq \cdots \leq ^p\lambda _j  \leq \cdots$ les valeurs propres   parabo\-liques rang\'ees par ordre croissant  avec multiplicit\'e.

\
 
  \noindent{\bf{Conjecture.}} Soit $S$ une surface hyperbolique d'aire finie. Alors on a    $\displaystyle ^p\lambda _{2g-2} >\frac 14$.
  
  \
  
Cette conjecture semble difficile \`a attaquer en utilisant uniquement les m\'ethodes topo\-lo\-giques de l'article. Il faudra probablement une meilleure compr\'ehension  du comportement pr\`es des pointes
de l'ensemble nodal d'une fonction somme de fonctions propres paraboliques. Une premi\`ere \'etape serait d'\'etudier
  si, comme celui d'une fonction $\lambda$-propre avec $\displaystyle \lambda \leq \frac 14$, l'ensemble nodal d'une fonction dans $\mathcal E$ est la r\'eunion d'un graphe dans la surface compl\'et\'ee $\bar S$ et d'un ensemble fini.  

\

Dans \cite{O}, on  proposait comme approche du probl\`eme de  la minoration de $\lambda _{2g-2}$, r\'esolu dans le Th\'eor\`eme \ref{2g-2}, 
la question interm\'ediaire suivante :  

\

  \noindent{\bf{Question 1.}} L'ensemble nodal d'une somme finie de fonctions propres non constantes est-il
un graphe incompressible (plus un ensemble discret)?  

\

Si la r\'eponse avait \'et\'e oui, on aurait pu reproduire mot-\`a-mot la d\'emonstration de \cite [Prop. 2]{O} pour obtenir le th\'eor\`eme 1. 
Mais c'est une question d'analyse apparemment difficile et
 qui n'est pas  purement locale puisque dans le rev\^etement universel du disque, il est facile de construire des contre-exemples.  Une autre question d'analyse,   locale cette fois est 
la suivante.

\

 \noindent{\bf{Question 2.}} Soit $S$ une surface munie d'une m\'etrique $C^\infty$. L'ensemble nodal d'une somme de fonctions  propres est-il la r\'eunion d'un graphe et d'un ensemble discret?

\

Essayons d'adapter l'argument de Cheng qui d\'ecrit l'ensemble nodal d'une fonction propre sur une surface. On remarque que si $\displaystyle f=\sum _1^kf_i$ est somme de fonctions $\lambda _i$-propres $f_i$, alors $f$ est la solution d'un op\'erateur elliptique, \`a savoir l'op\'erateur $ \displaystyle L= (\Delta +\lambda _1)\circ \cdots \circ (\Delta  +\lambda_k)$. L'ordre d'annulation de $f$  au point $p$ est fini par un th\'eor\`eme de Carleman \cite {Ca}, \cite{A}. Un th\'eor\`eme de Bers \cite {Be}  pr\'ecise   alors la partie principale $P$ du d\'eveloppement de Taylor au point $p$. Cette partie principale est   solution de l'equation  $L_0 (P)=0$, o\`u $L_0$ est la partie principale de l'op\'erateur $L$ au point $p$ : dans une carte conforme autour de $p$, on a  $L_0=\Delta _{eu}^k$ 
o\`u
$\Delta _{eu}$ est le Laplacien euclidien. Donc la partie principale $P$ de  $f$
 est un polyn\^ome particulier, dit {\it polyn\^ome polyharmonique}. Mais ce contr\^ole de la partie principale de $f$ au point  $p$ n'est pas suffisant pour pouvoir conclure  que  l'ensemble nodal de $f$ est isotope \`a celui de $P$ comme dans la d\'emonstration du th\'eor\`eme de Cheng.   En effet, pour construire l'isotopie (cf. \cite {K}), 
 $P$ doit v\'erifier une propri\'et\'e de transversalit\'e  sur son ensemble nodal : son gradient ne doit pas s'y annuler   (sauf peut-\^etre au point $p$). Or, cette propri\'et\'e qui est v\'erifi\'ee par un polyn\^ ome harmonique (la  partie principale dans le cas d'un seule fonction propre) ne l'est plus n\'ecessairement pour les polyn\^omes polyharmoniques.  
 C'est pour cette raison que  le th\'eor\`eme \ref{2g-2} contient l'hypoth\`ese que la m\'etrique \'etudi\'ee est analytique.
Toutefois, une r\'eponse positive \`a la Question 2 permettrait de  remplacer ``analytique'' par ``$C^\infty$'' dans le th\'eor\`eme \ref{2g-2}.


\begin{thebibliography}{EMOT}
 
\bibitem[A]{A}N. Aronszajn, 
{\em Sur l'unicit\'e du prolongement des solutions des \'equations aux d\'eriv\'ees partielles elliptiques du second ordre},  
C. R. Acad. Sci. Paris 242 (1956), 723--725. 
 
 
 \bibitem[Be]{Be} L. Bers,    {\em  Local behaviour of solutions of general linear elliptic equations}, 
 Comm. Pure Appl. Math. 8 (1955), 475--504.

\bibitem[Br]{Br} R. Brooks, 
{\em  The fundamental group and the spectrum of the Laplacian},
Comment. Math. Helv. 56 (1981), no. 4, 581--598. 


 
\bibitem[Bu1]{Bu1} P. Buser,  
{\em Riemannsche Fl\"achen mit Eigenwerten in $(0,$ $1/4)$}, 
Comment. Math. Helv. 52 (1977), 25--34. 
 
 
\bibitem[Bu2]{Bu2} P. Buser, {\em Geometry and spectra of compact Riemann surfaces},     Birkha\"user Verlag, Basel, Berlin, Boston MA,  (1992). 
 
 \bibitem[Ca]{Ca} T. Carleman,   {\em Sur les syst\`emes lin«eaires aux 
 d\'eriv\'ees partielles du premier ordre \`a
deux variables}, C. R. Acad. Sc. Paris 197 (1933), 471--474.
  
  
 \bibitem[Cha] {Cha} I. Chavel,  {\em  Eigenvalues In Riemannian Geometry}, Academic Press, Orlando,
  (1984).
 
 
 \bibitem[Che]{Che} S.-Y. Cheng,  {\em Eigenfunctions and nodal sets}, 
Comment. Math. Helvetici  51 (1976),  43--55.
 
 

\bibitem[Hu] {Hu} M. N. Huxley, 
 {\em Cheeger's Inequality with a boundary term},
  Commentarii Mathematici Helvetici  58  (1983), 347--354.
 
 

\bibitem[I] {I} H. Iwaniec, 
{\em Introduction to the spectral theory of automorphic forms},
Biblioteca de la Revista Matem\'atica Iberoamericana 
Madrid (1995).


\bibitem[K] {K} T-C. Kuo,
{\em On $C^0$-sufficiency of jets of potential functions
}, Topology    8 (1969),  167--171.
 
 

 
 

 \bibitem[Lo] {Lo}  
 S. Lojasiewicz,  {\em  Sur le probl\`eme de la division}, 
Studia Mathematica    28   (1959),   87--136.

  
  \bibitem [O]{O}    J.-P. Otal,
{\em   Three topological properties of small eigenfunctions on
hyperbolic surfaces},
Geometry and Dynamics  of Groups and Spaces,  in Memory of Alexander Reznikov,   Birkha\"user,  Progress in Mathematics  265, (2008).
     


 \bibitem[PS] {PS} 
 R. Phillips et P. Sarnak,  {\em  On cusps forms for cofinite subgroups of $\rm{PSL(2},\RR)$}, 
Invent. Math. 80 (1985), 339--364.


 
 
  \bibitem [R1]{R1} 
B. Randol, 
{\em Small eigenvalues of the Laplace operator on compact Riemann surfaces},
  Bull. Am. Math. Soc. 80  (1974) 996-1000. 

  \bibitem [R2]{R2} 
B.  Randol,   
{\em A remark on $\lambda \sb{2g-2}$},
 Proc. Am. Math. Soc. 108, No.4, (1990)  1081-1083 (1990). 
   
  
  \bibitem[Sc]{Sc}
P. Schmutz,  
{\em Small eigenvalues on Riemann surfaces of genus $2$},
Invent. Math. 106 (1991),  121--138.
  
  
  
  \bibitem [SWP]{SWP}
R. Schoen, S. Wolpert, S.T. Yau, 
{\em Geometric bounds on the low eigenvalues of a compact surface,} Geometry of the Laplace operator (Proc. Sympos. Pure Math., Univ. Hawaii, Honolulu, Hawaii, 1979),  pp. 279--285, 
Proc. Sympos. Pure Math., XXXVI, Amer. Math. Soc., Providence, R.I., 1980. 
  
  

\bibitem [Se]{Se}   B. S\'evennec,
 {\em  Multiplicity of the second  
Schr\"odinger eigenvalue on closed surfaces},
Math. Ann.
  324    (2002),  195--211.
   
  
\bibitem [Z]{Z}    P. Zograf, 
{\em   Small eigenvalues of automorphic Laplacians in spaces of cusp forms},
 Soviet Math. Dokl.  27 (1983),   420--422.
     



  
\end{thebibliography}
\end{document}